\documentclass[journal,onecolumn,11pt]{IEEEtran}

\input packages-headers
\RequirePackage[top=3.5cm,bottom=3.5cm,left=3.5cm,right=3.5cm,nohead]{geometry}

\usepackage[american]{circuitikz}

\newcommand{\captionfonts}{\small}

\begin{document}

\pagestyle{plain}

\title{\huge Analyzing Oscillators using Describing Functions}

\author
{
Tianshi Wang\\
\vspace{0.5em}
{\small
Department of EECS, University of California, Berkeley, CA, USA \\
Email: \texttt{tianshi@berkeley.edu}}
}

\maketitle
\thispagestyle{plain}

\begin{abstract}
In this manuscript, we discuss the use of describing functions as a systematic
approach to the analysis and design of oscillators.
Describing functions are traditionally used to study the stability of nonlinear
control systems, and have been adapted for analyzing LC oscillators.
We show that they can be applied to other categories of oscillators too,
including relaxation and ring oscillators.
With the help of several examples of oscillators from various physical domains,
we illustrate the techniques involved, and also demonstrate the effectiveness
and limitations of describing functions for oscillator analysis.
\end{abstract}

\section{Introduction} \seclabel{intro}

Oscillators can be found in virtually every area of science and technology.
They are widely used in radio communication \cite{rappaport1996wireless}, clock
generation \cite{best2007PLL}, and even logic computation in both Boolean
\cite{WaRoUCNC2014PHLOGON,WaRoDAC2015MAPPforPHLOGON,WaRoOscIsing2017,Wa2017arXivMetronomes}
and non-Boolean
\cite{Hoppensteadt2000synchronization,maffezzoni2015oscArray,Porod2015physical}
paradigms.
The design and use of oscillators are not limited to electronic ones.
Integrated MEMS oscillators are designed with an aim of replacing traditional
quartz crystal ones as on-chip frequency references
\cite{miri2010design,nguyen1999integrated,nguyen1995micromechanical};
spin torque oscillators are studied for potential RF applications
\cite{dixit2012spintronic,villard2010ghz,finocchio2013nanoscale};
chemical reaction oscillators work as clocks in synthetic biology
\cite{danino2010synchronized,rosenblum2003synchronization};
lasers are optical oscillators with a wide range of applications.
The growing importance of oscillators in circuits and systems calls for a
systematic design methodology for oscillators which is independent from
physical domains.
However, such an approach is still lacking.
While structured design methodologies are available for oscillators with
negative-feedback amplifiers \cite{verhoeven2010structured,van2001structured},
different types of oscillators are still analyzed with vastly different
methodologies.
The fundamental design specifications, such as signal waveforms and frequency
tuning range, \etc, are all calculated using separate models;
the interconnections between them can be easily overlooked during design,
causing the iterations and time to increase.
Partly as a result of this lack of structured study, researchers in domains
other than electronics are often not aware of the similarity between their
oscillators and some of the existing electronic ones.
So instead of applying the existing models developed in electronics, they often
rely heavily on experiments to guide their design, leading to higher design
costs as well as less optimal results.

In this manuscript, we use describing function analysis to study oscillators in
a systematic way.
Describing function analysis has been practically applied to nonlinear control
system design for many decades \cite{vander1968multiple}.
It is a general approach for analyzing the stability as well as predicting
limit cycle properties such as frequency and amplitude of nonlinear systems.
It generates graphical results in the complex plane that can provide system
designers with convenient visualization for analyzing oscillation and exploring
design margins.
It is an immediate and natural candidate for studying oscillators in a
systematic way.
And in fact, it has been successfully applied to oscillator design previously
\cite{bank2006harmonic,vidal2001describing}.
However, these works focus mainly on negative-resistance-feedback LC
oscillators in the electronic domain.
Ring oscillators are never part of the scope.
Neither are relaxation oscillators, especially when the amplifier inside has
non-negligible delays.
Overall, the use of describing function in the context of general
domain-independent oscillator design still warrants discussion; examples from
oscillators in different physical domains are highly needed.

We start the manuscript by providing the readers with essential background
knowledge of describing function analysis in the context of nonlinear control
theory in \secref{background}.
Then we use the analysis to study three common categories of oscillators:
negative-resistance-feedback LC oscillators, relaxation oscillators and ring
oscillators in \secref{analysis}.
Examples from electronics, biology and neuroscience are used as demonstration.
In the meanwhile, we also address the limitations of describing function method
and provide thoughts on potential improvements on it for oscillator design.
Furthermore, in \secref{design}, we explore the use of describing function in
the design of oscillators.
In particular, we redesign a relaxation oscillator to make its response more like
a harmonic oscillator in Section \secref{relaxationdesign}.
We show that despite differences in their physical implementations, relaxation
and harmonic LC oscillators can be treated as the same type of oscillators from
describing function analysis' point of view, and they share the same
mathematical equations.
Conclusions and further research plans are provided in Section
\secref{conclusions}.

\section{Background} \seclabel{background}

In this section we provide background knowledge of describing function
analysis in the context of nonlinear control theory.

\subsection{What are describing functions?}

Given a general nonlinear algebraic function $y = f(x)$, without being
concerned with its physical interpretation or application, we can define and
derive its describing function in the following manner:

\begin{enumerate}
    \item set $x(A, t) = A \cdot \sin(\omega t)$ as the input to $f(x)$,
          $x(A, t)$ is periodical in $t$ with period $2\pi/\omega$;
    \item get output: $y(A, t) = f(x(A, t)) = f(A \cdot \sin(\omega t))$,
          $y(A, t)$ is also periodical in $t$ with period $2\pi/\omega$;
    \item Fourier transform $y(A, t)$, get Fourier coefficient $Y(A)$ at
          fundamental frequency $\omega/2\pi$;
    \item define the describing function of $y = f(x)$ as $N(A) = Y(A)/A$.
\end{enumerate}

Mathematically, we can write $N(A)$ as
\begin{equation} \label{eq:NA}
  N(A) = \frac{1}{A} (a_1 + j\cdot b_1),
\end{equation}
where
\begin{eqnarray}
  a_1 &=& +\frac{1}{\pi}\int_0^{2\pi} f(A \cdot \sin(\omega t)) \sin(\omega t) d\omega t \label{eq:a_1},\\
  b_1 &=& -\frac{1}{\pi}\int_0^{2\pi} f(A \cdot \sin(\omega t)) \cos(\omega t) d\omega t \label{eq:b_1}.
\end{eqnarray}

Note that the choice of $\omega$ doesn't affect the value of $N(A)$ as $f(x)$
is a memoryless algebraic function.
This can be explained in (\ref{eq:a_1}) and (\ref{eq:b_1}). 
When we consider $\omega t$ as a single variable, the values of the integrals
don't rely on the choice of the notation of this variable.

Graphically, the definition of $N(A)$ can be illustrated as in \figref{DFdefinition}.
With sinusoidal input $x(A, t) = A\cdot \sin(\omega t)$, $Y(A)$ can be viewed
as the amplitude of the sinusoidal wave at the fundamental frequency of the
output.
In this way, the describing function can be intuitively thought of as the
frequency domain linear approximation of the nonlinear algebraic function that
transfers an input phasor (sinusoidal wave) with amplitude $A$ to an output
phasor $Y(A) = A\cdot N(A)$.

\begin{figure}[htbp]
\centering{
	\epsfig{file=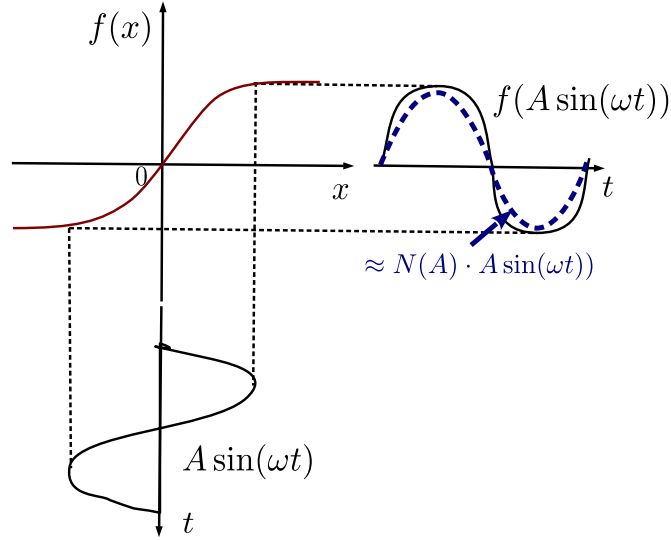,width=0.60\textwidth}
}
\caption{Plot illustrating the definition of describing function. \figlabel{DFdefinition}}
\end{figure}

\subsection{How to calculate describing functions?}

The definition of describing function as in (\ref{eq:NA}), (\ref{eq:a_1}) and
(\ref{eq:b_1}) provides a direct method for its calculation.
Analytically calculated describing functions for some nonlinear functions
common in control systems are listed in Table \ref{tab:DFs}.
Take the ideal saturation characteristic function as shown in Table
\ref{tab:DFs} No. 1 as an example:
\begin{equation}\label{eq:fsat}
    f(x) =
    \begin{cases}
        -K\cdot a &\mbox{if }  x \leq -a \\
         K\cdot x &\mbox{if }  -a < x < a \\
         K\cdot a &\mbox{if }  x \ge +a.
    \end{cases}
\end{equation}

With input $x(A, t) = A\cdot \sin(\omega t)$, output $y(A, t)$ is
\begin{equation}
    y(A, t) =
    \begin{cases}
        Ka &\mbox{if }   2k\pi + \Phi < \omega t < (2k+1)\pi - \Phi \\
        -Ka &\mbox{if }  (2k+1)\pi + \Phi < \omega t < (2k+2)\pi - \Phi \\
        KA\sin(\omega t) &\mbox{everywhere else,}
    \end{cases}
\end{equation}
where $k\in \mathbf{Z}$ and $\Phi = \arcsin(a/A)$, assuming $A \ge a$.

Now we calculate $N(A)$ analytically.
\begin{eqnarray}
  a_1 &=& \frac{1}{\pi}\int_0^{2\pi} f(A \cdot \sin(\omega t)) \sin(\omega t) d\omega t \\
      &=& \frac{4}{\pi}\int_0^{\pi/2} f(A \cdot \sin(\omega t)) \sin(\omega t) d\omega t \\
      &=& \frac{4}{\pi}\left[\int_0^{\Phi} KA\sin^2(\omega t)  d\omega t
          + \int_{\Phi}^{\pi/2} Ka \sin(\omega t) d\omega t \right] \\
      &=& \frac{2KA}{\pi}\left[ \arcsin (\frac{a}{A}) +
          \frac{a}{A}\sqrt{1 - (\frac{a}{A})^2}\right].
\end{eqnarray}

Similarly, we have
\small{
\begin{eqnarray}
  b_1 &=& -\frac{1}{\pi}\int_0^{2\pi} f(A \cdot \sin(\omega t)) \cos(\omega t) d\omega t \\
      &=& -\frac{1}{\pi}\int_0^{\pi/2} f(A \cdot \sin(\omega t)) \cos(\omega t) d\omega t 
          -\frac{1}{\pi}\int_{\pi/2}^{\pi} f(A \cdot \sin(\omega t)) \cos(\omega t) d\omega t \nonumber\\
       && -\frac{1}{\pi}\int_{\pi}^{3\pi/2} f(A \cdot \sin(\omega t)) \cos(\omega t) d\omega t 
          -\frac{1}{\pi}\int_{3\pi/2}^{2\pi} f(A \cdot \sin(\omega t)) \cos(\omega t) d\omega t.
\end{eqnarray}
}

Because
\small{
\begin{eqnarray}
  && \frac{1}{\pi}\int_{\pi/2}^{\pi} f(A \cdot \sin(\omega t)) \cos(\omega t) d\omega t \nonumber \\
  &=& \frac{1}{\pi}\int_{\pi/2}^0 f(A \cdot \sin(\pi-\omega t)) (-\cos(\pi-\omega t)) (-d(\pi-\omega t)) \\
  &=& -\frac{1}{\pi}\int_0^{\pi/2} f(A \cdot \sin(\theta)) \cos(\theta) d\theta,
\end{eqnarray}
}
where $\theta = \pi-\omega t$.

We have
\begin{equation}
    \frac{1}{\pi}\int_0^{\pi/2} f(A \cdot \sin(\omega t)) \cos(\omega t) d\omega t 
	 +\frac{1}{\pi}\int_{\pi/2}^{\pi} f(A \cdot \sin(\omega t)) \cos(\omega t) d\omega t = 0,
\end{equation}
and similarly
\begin{equation}
	\frac{1}{\pi}\int_{\pi}^{3\pi/2} f(A \cdot \sin(\omega t)) \cos(\omega t) d\omega t 
	+\frac{1}{\pi}\int_{3\pi/2}^{2\pi} f(A \cdot \sin(\omega t)) \cos(\omega t) d\omega t = 0.
\end{equation}

Therefore
\begin{equation}
 \implies b_1 = 0.
\end{equation}

Finally the analytical expression of the describing function of ideal saturation
function (\ref{eq:fsat}) can be calculated as

\begin{eqnarray}
  N(A) &=& \frac{1}{A} (a_1 + j\cdot b_1) \nonumber \\
  &=& \frac{2K}{\pi}\left[ \arcsin (\frac{a}{A}) +
          \frac{a}{A}\sqrt{1 - (\frac{a}{A})^2}\right], (A\ge a)
\end{eqnarray}

Describing functions shown in Table \ref{tab:DFs} are calculated
analytically in the same manner.

However, we may not always be able to express analytical integrals as in
\ref{eq:a_1} and \ref{eq:b_1} explicitly.
Alternative to the analytical calculation procedures, describing functions
can also be calculated numerically just by sweeping $A$ of the input, generating
outputs through interpolation on tabulated data of $f(x)$, and performing
Fourier transformation on the outputs to acquire the Fourier coefficients at
the fundamental frequency.

\begin{table}[htb]
	\begin{center}
		\begin{tabular}{ |m{0.06\textwidth}|m{0.2\textwidth}|m{0.7\textwidth}|}
		\hline
		{\bf No.} & {\bf Non-linearities} & {\bf Describing Functions} \\ \hline
		\hline
		{\bf 1} & 
		 \centering{
		   \epsfig{file=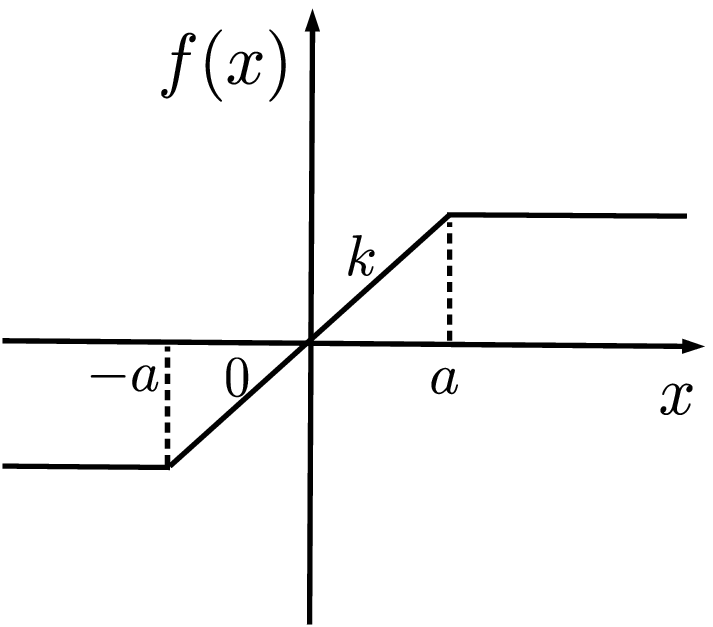,width=0.20\textwidth}
		 }
		& $N(A) = \frac{2K}{\pi}\left[ \arcsin (\frac{a}{A}) + \frac{a}{A}\sqrt{1 - (\frac{a}{A})^2}\right]$, $(A\ge a)$
		\\ \hline
		{\bf 2} & 
		 \centering{
		   \epsfig{file=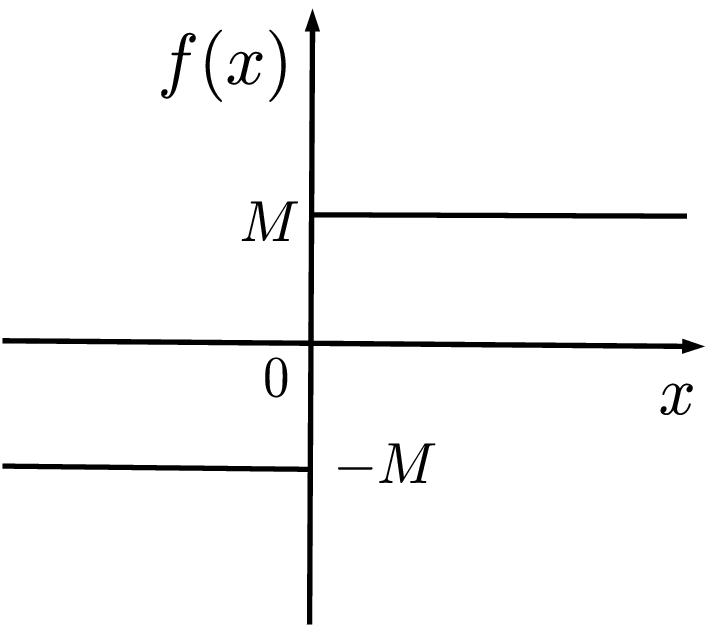,width=0.20\textwidth}
		 }
		& $N(A) = \frac{4M}{\pi A}$
		\\ \hline
		{\bf 3} & 
		 \centering{
		   \epsfig{file=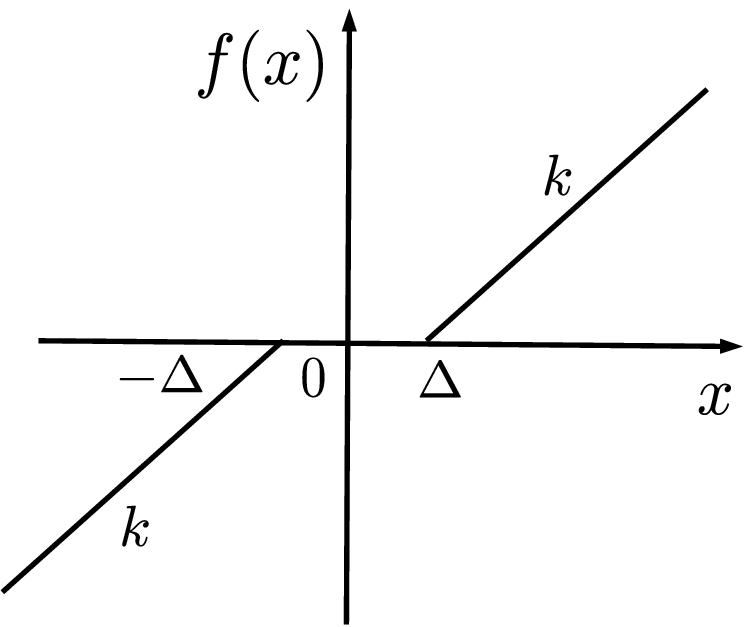,width=0.20\textwidth}
		 }
		& $N(A) = \frac{2K}{\pi}\left[ \frac{\pi}{2} - \arcsin(\frac{\Delta}{A})
		  - \frac{\Delta}{A}\sqrt{1 - (\frac{\Delta}{A})^2}\right]$, $(A\ge \Delta)$
		\\ \hline
		{\bf 4} & 
		 \centering{
		   \epsfig{file=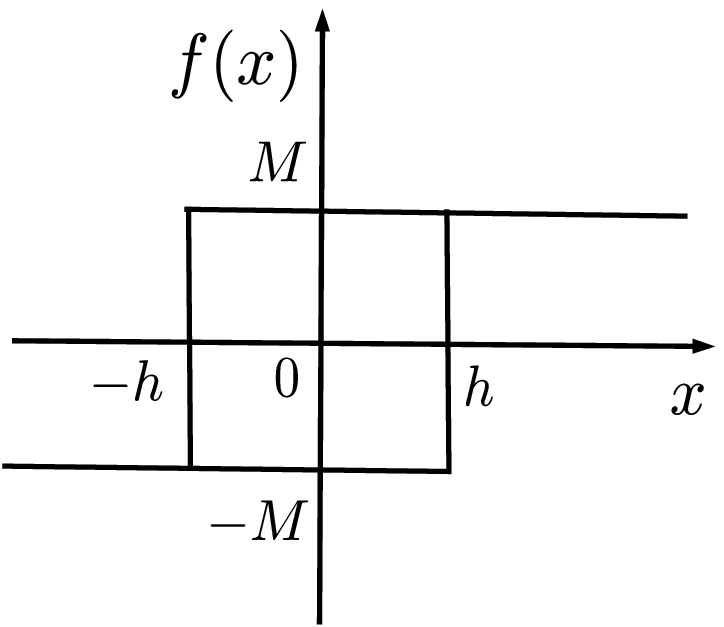,width=0.20\textwidth}
		 }
		& $N(A) = \frac{4M}{\pi A} \sqrt{1-(\frac{h}{A})^2} - j \frac{4Mh}{\pi A^2}\cdot$, $(A \ge h)$
		\\ \hline
		\end{tabular}
	\end{center}
  \caption{Describing Functions for nonlinearities that are common in control theory. \label{tab:DFs}}
\end{table}

\subsection{How to use describing functions?}

Describing function $N(A)$ of a nonlinear function in the loose sense transfers 
a sinusoidal wave with amplitude $A$ to another sinusoidal wave with amplitude
$|N(A)\cdot A|$ and phase $\angle (N(A))$.
In this way, describing function $N(A)$ can be viewed as the frequency domain
approximation of $f(x)$ that maps input phasor $A$ to output phasor $N(A)\cdot A$.
Note that when $f(X)$ is linear, analytical calculation shows that its
describing function is constant regardless of input amplitude $A$.
This is consistent with the frequency domain expression of one-port linear systems.

With describing functions acting as approximations to nonlinear functions, we can
analyse general nonlinear systems in frequency domain and apply stability criteria
that are available in linear control theory.


\tikzstyle{block} = [draw, fill=blue!20, rectangle, 
    minimum height=3em, minimum width=6em]
\tikzstyle{sum} = [draw, fill=blue!20, circle, node distance=1cm]
\tikzstyle{input} = [coordinate]
\tikzstyle{output} = [coordinate]
\tikzstyle{pinstyle} = [pin edge={to-,thin,black}]

\begin{figure}[htbp]
\centering{
\begin{tikzpicture}[auto, node distance=2cm,>=latex']
    \node [input, name=input] {};
    \node [sum, right of=input] (sum) {};
    \node [block, right of=sum] (linear) {$G(s)$};
    \node [block, right of=linear, node distance=3cm] (nonlinear) {$N(A)$};
    \draw [->] (linear) -- node[name=x] {$X(s)$} (nonlinear);
    \node [output, right of=nonlinear] (output) {};
    \coordinate [below of=x] (tmp);

    \draw [draw,->] (input) -- node {$U(s)$} (sum);
    \draw [->] (sum) -- node {$E(s)$} (linear);
    \draw [->] (nonlinear) -- node [name=y] {$Y(s)$}(output);
    \draw [->] (y) |- (tmp) -| node[pos=0.99] {$-$} 
        node [near end] {} (sum);
\end{tikzpicture}
}
\caption{{\captionfonts A general nonlinear feedback system with separated 
     one-port linear and nonlinear blocks.
	 \figlabel{nonlinear_feedback}}}
\end{figure}
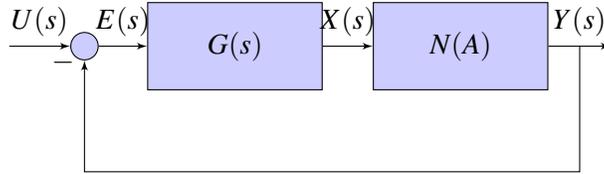

Consider a nonlinear feedback system whose nonlinear and linear parts can be
separated, as is show in \figref{nonlinear_feedback}.
Using describing function, we write its transfer function as:
\begin{equation}\label{eq:H}
	H(s) = \frac{G(s)\cdot N(A)}{1 + G(s)\cdot N(A)}
\end{equation}

This is a transfer function with a parameter $A$. 
At any given $A$, we can assess the stability of this system using Nyquist
stability criterion.
However, this lead to a common confusion of people new to describing function
analysis: if $N(A)$ is only a valid approximation when the input and output of
$f(x)$ are both sinusoidal with certain amplitudes, what does it mean when we
say a system is stable or unstable?

Describing function analysis first assumes the system is already oscillating
and periodical waveforms present everywhere in the system are all sinusoidal.
With a given $A$, we further assumes the oscillating waveform at the input
of the nonlinear block has amplitude of $A$.
Under all these assumptions, the transfer function in (\ref{eq:H}) holds well.
Nyquist stability criterion is then applied to check whether the system is
stable.
If it is strictly stable (all closed-loop poles are in left-open plane) or 
strictly unstable, it indicates the assumptions we made are violated.
Only when the system is marginally stable can there be chance for the
system to actually oscillate at the assumed amplitude.
The necessary condition for the system to be marginally stable is
\begin{equation}\label{eq:marginalstable}
	-\frac{1}{N(A)} = G(s)
\end{equation}

The condition (\ref{eq:marginalstable}) is satisfied when the plot of
$-\frac{1}{N(A)}$ on complex plane graphically intersects with the Nyquist plot
of $G(j\omega)$.
And the $\omega$ and $A$ with respect to the intersection provide predictions
to the oscillation's frequency and amplitude (at the input of nonlinear block)
respectively. 
As an example, in \figref{nyquist_stability} there are two intersections ---
$A_1$ and $A_2$, indicating two possible oscillating operating points.
Compared with other nonlinear system analyses like Harmonic Balance,
describing function analysis not only locates possible oscillating operating
points, but also provides information for the stability of oscillation.
In \figref{nyquist_stability}, if we assume that the system is oscillating at
operating point $A_1$.
If there is some perturbation to the input of nonlinear block, $A$ may change
by a small amount. 
When $A$ increases, $A_1$ shifts to $B_1$; the resulting closed-loop system
becomes stable, causing $A$ to decay and the operating point to move back to
$A_1$.
When $A$ decreases slightly to reach $C_1$ where the closed-loop system is
unstable, $A$ increases back to its value at the operating point $A_1$.
These observations indicate that the oscillation happening at intersection $A_1$
is stable.
On the contrary, similar analysis predicts that oscillation at the other
intersection $A_2$ is unstable, thus cannot be observed in the actual system.

\begin{figure}
 \centering{
   \epsfig{file=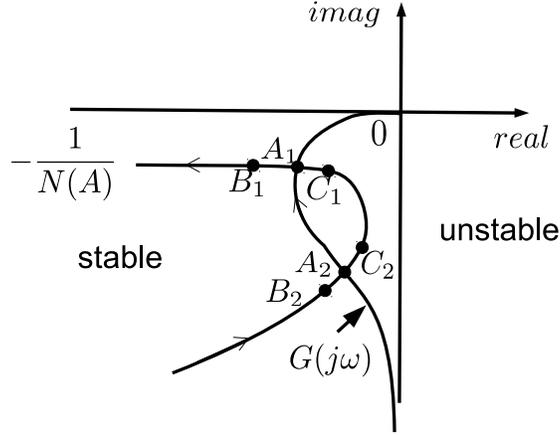,width=0.50\textwidth}
 }
 \caption{{\captionfonts Nyquist plot of $G(j\omega)$ intersected with $-\frac{1}{N(A)}$.
	 Arrows on $G(j\omega)$ and $-\frac{1}{N(A)}$ indicate the increasing
     directions of $\omega$ and $A$ respectively.
     When $-\frac{1}{N(A)}$ lands in area left of the Nyquist plot of $G(j\omega)$, the
     closed-loop system is stable.
     Two intersections are $A_1$ and $A_2$.
	 \figlabel{nyquist_stability}}}
\end{figure}

\section{Oscillator Analysis with Describing Functions} \seclabel{analysis}

There are three main categories of oscillators: LC oscillators, ring
oscillators and relaxation oscillators.
In this section, we apply describing function analysis to all of them.
In each category, we select typical examples, identify and separate linear and
nonlinear parts of the systems, then apply describing function analysis to
them and study its effectiveness and limitations.

\subsection{LC Oscillators} \seclabel{LC}

Applying describing function analysis to negative-resistance-feedback LC
oscillators is a natural choice, and has been studied before
\cite{bank2006harmonic,vidal2001describing}.
In this section, we briefly go through the procedures and analyze an LC
oscillator with a less common series RLC configuration, just to show the
generality of the describing function analysis.

It is worth noting that a wide range of oscillators can be categorized as
LC-type oscillators, including spin torque oscillators, MEMS oscillators,
optical lasers, \etc{}
The same analysis can be applied to all of them.

A schematic diagram of a series LC oscillator is given in the left-most
plot in \figref{LC}. 
The negative resistor is usually nonlinear to ensure that the oscillation is
amplitude-stable.
To apply describing function analysis, we separate linear and nonlinear
parts of the system as shown in the middle plot of \figref{LC} and convert the
system into block diagram.
\begin{figure}[htbp]
\centering{
	\epsfig{file=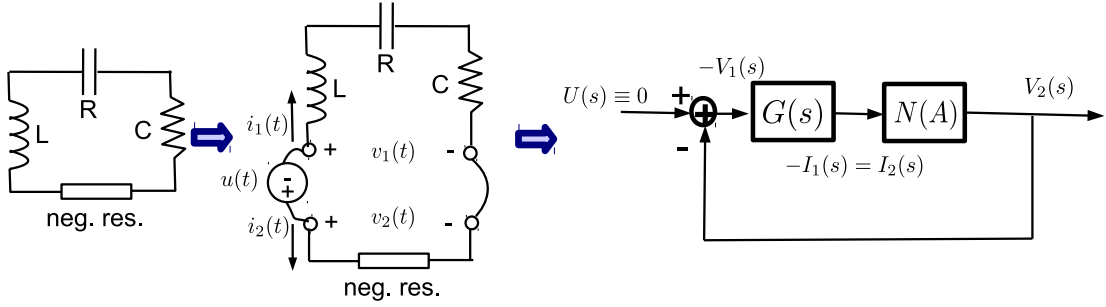,width=\textwidth}
}
\caption{Schematics and block diagram of a negative-resistance-feedback LC
     oscillator.
	 The middle diagram shows how the series R-L-C-negative-resistor
	 system is analyzed as a connection of nonlinear and linear parts.
	 \figlabel{LC}}
\end{figure}

In the block diagram in \figref{LC}:
\begin{equation} \label{eq:LCG}
	G(s) = \frac{1}{R+\frac{1}{sC} + sL} = \frac{sC}{LC\cdot s^2 + RC\cdot s + 1}.
\end{equation}

We can assume the negative nonlinear resistor has the following characteristic:
\begin{equation} \label{eq:LCf}
	v_2 = f(i_2) = -V_{max} \cdot \tanh (\frac{R_{max}}{V_{max}} \cdot i_2).
\end{equation}

The describing function of (\ref{eq:LCf}) can be acquired numerically by
sweeping $A$ in the sinusoidal input $i_2 = A\cdot \sin(\omega t)$. 
However, just to show a flavour of how describing functions can be
calculated by hand in general, here we derive the describing function of
(\ref{eq:LCf}) using Taylor expansion.
\begin{equation}
	N(A) = \frac{1}{A\pi} \int_0^{2\pi} -V_{max} \cdot \tanh (\frac{R_{max}}{V_{max}} \cdot
	A\sin(\omega t)) \cdot \sin(\omega t) d(\omega t).
\end{equation}

Because $\tanh(x) = x - \frac{1}{3}x^3 + \frac{2}{15}x^5 + o(x^5)$, we have
\begin{eqnarray}
	N(A) = && \frac{1}{A\pi} \int_0^{2\pi}
	-V_{max} \cdot
	\frac{AR_{max}}{V_{max}} \cdot \sin^2(\omega t)
	- \frac{1}{3} \frac{A^3 R^3_{max}}{V_{max}} \cdot \sin^4(\omega t) \nonumber \\
	&& +\frac{2}{15} \frac{A^5 R^5_{max}}{V_{max}} \cdot \sin^6(\omega t)
	+o(\sin^6(\omega t))
	d(\omega t). \label{eq:LCNcomplicated}
\end{eqnarray}

Using
\begin{eqnarray}
	\int\sin^2(x) dx &=& \frac{1}{2} (1-\cos(2x)), \\
	\int\sin^4(x) dx &=& \frac{1}{8} (3-4\cos(2x)+\cos(4x)), \\
	\int\sin^6(x) dx &=& \frac{1}{32} (10-15\cos(2x)+6\cos(4x)-\cos(6x)),
\end{eqnarray}
(\ref{eq:LCNcomplicated}) is reduced to
\begin{equation} \label{eq:LCN}
	N(A) = 
	-R_{max}
	+\frac{R^3_{max}}{4V^2_{max}} \cdot A^2
	-\frac{R^5_{max}}{12V^4_{max}} \cdot A^4
	+ o(A^4).
\end{equation}

After characterizing the nonlinear block using describing function, we derive the 
transfer function of the system and the condition for it to oscillate.
The transfer function is
\begin{equation}
	H(s) = \frac{V_2(s)}{U(s)} = \frac{G(s)\cdot N(A)}{1+G(s)\cdot N(A)},
\end{equation}
where $G(s)$ and $N(A)$ are given in (\ref{eq:LCG}) and (\ref{eq:LCN}) respectively.
For the system to oscillate, $-1/N(A) = G(s)$ has to be satisfied.
Given specific parameters, we will then be able to location the oscillation
operating point of the system and predict the oscillation frequency as well as
amplitudes of the waveforms.

\subsection{Ring Oscillators} \seclabel{ring}

Ring oscillators are another type of oscillators that normally consist of odd
numbers of inverters.
The functionality of this type of oscillators relies on the delay each stage of
inversion provides.
Therefore, it is less obvious how memoryless nonlinearity can be separated from
the system and how describing function analysis can be used. 
However, with some assumptions, describing function analysis can still be
applied in this scenario.

\begin{figure}[htbp]
\centering{
	\epsfig{file=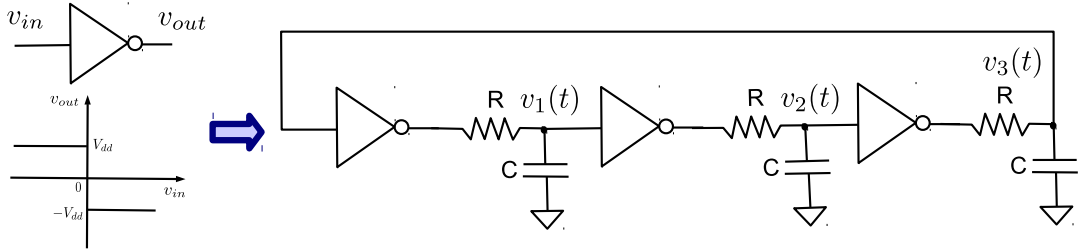,width=\textwidth}
}
\caption{ 3-stage ring oscillator made of ideal inverters.
    \figlabel{ringosc}}
\end{figure}

\figref{ringosc} shows a three-stage ring oscillator made of inverters with an
ideal inverse relay characteristics.
It is common to assume that the three stages are identical, in both the
inverter and the RC delay.
The analytical solution for the oscillation frequency and amplitude of this
system can be derived \cite{SrRoTCAS2007}:
\begin{eqnarray}
	A^* &=& \frac{\sqrt{5}-1}{2} \approx 0.618 \label{eq:ringAstar},\\
	T^* &=& 6\ln(\frac{1 + \sqrt{5}}{2}) \cdot \tau \approx 2.89\tau, \label{eq:ringTstar}
\end{eqnarray}
where $\tau = RC$. $V_{dd}$ is normalized to 1.

Now we apply describing function analysis to the system.
We treat the first inverter (leftmost in \figref{ringosc}) as a memoryless
nonlinearity.
According to Table \ref{tab:DFs}, we know its describing function is
\begin{equation} \label{eq:ringN}
	N(A) = -\frac{4V_{dd}}{\pi A}
\end{equation}

The RC circuit can be easily modeled as a linear block.
Since the stages are identical, the signals after the three stages are of the
same waveforms with phases separated by $2\pi/3$.
The waveform of $v_3(t)$ can be thought of as a time-delayed version of
$v_1(t)$.
Therefore, we can model the rest of the two stages as a pure time delay, which
is a linear subsystem.
This technique is illustrated in \figref{ringosc2}.
After these procedures, we can draw the block diagram of this ring oscillator
as in \figref{ringosc3}.
\begin{figure}[htbp]
\centering{
	\epsfig{file=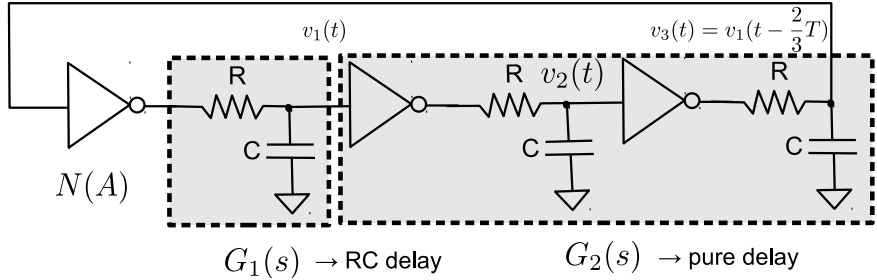,width=0.8\textwidth}
}
\caption{ Modelling the ring oscillator in \figref{ringosc} by identifying nonlinear
    and linear blocks in the system.
    \figlabel{ringosc2}}
\end{figure}

\begin{figure}[htbp]
\centering{
	\epsfig{file=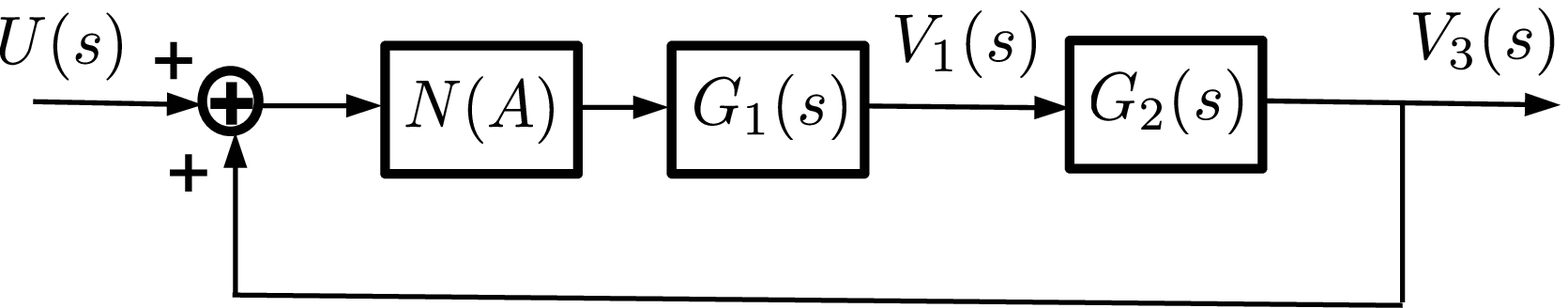,width=0.7\textwidth}
}
\caption{ Block diagram of the system in \figref{ringosc}.
    \figlabel{ringosc3}}
\end{figure}

The transfer function derived from the block diagram in \figref{ringosc3} is
\begin{equation}
	H(s) = \frac{V_3(s)}{U(s)} = \frac{G_2(s) \cdot G_1(s) \cdot N(A)}{1 - G_2(s) \cdot G_1(s) \cdot N(A)},
\end{equation}
where $N(A)$ is given in (\ref{eq:ringN}), $G_1(s)$ and $G_2(s)$ are
\begin{equation} \label{eq:ringG1G2}
	G_1(s) = \frac{1}{RC\cdot s + 1}, \text{~~}
	G_2(s) = e^{-\frac{2}{3}Ts}.
\end{equation}

For the system to oscillate, condition $1/N(A) = G_2(s)\cdot G_1(s)$ has to be
satisfied.
We can examine this condition by overlaying the Nyquist plot of linear system
$G_2(s)\cdot G_1(s)$ and the plot of $1/N(A)$ in the same complex plane in
\figref{ring_nyquist}.
\begin{figure}[htbp]
\centering{
	\epsfig{file=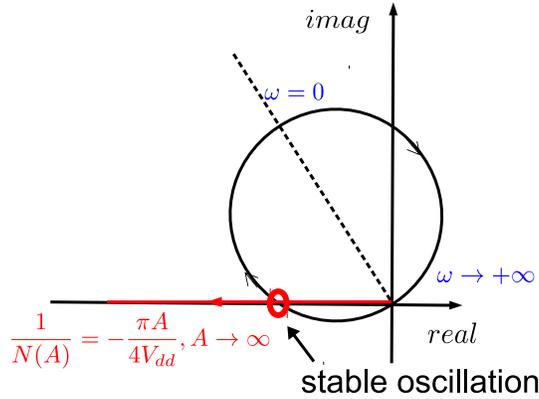,width=0.5\textwidth}
}
\caption{ Nyquist plot of $G_2(s)\cdot G_1(s)$ as in (\ref{eq:ringG1G2}) and
the plot of $1/N(A)$ as in (\ref{eq:ringN}) overlaid in the same complex plane.
The intersection indicates a stable oscillating operating point.
    \figlabel{ring_nyquist}}
\end{figure}

Mathematically, the oscillation intersection in \figref{ring_nyquist} implies
the following condition.
\begin{equation}
	\frac{1}{N(A)} = G_2(j\omega)\cdot G_1(j\omega) = \frac{1}{RC\cdot j\omega + 1} e^{-j\frac{4\pi}{3}}.
\end{equation}

From \figref{ring_nyquist} we observe that the intersection happens at $(-0.5,~
0)$, which indicates $1/N(A) = -1/2$.
\begin{equation} \label{eq:ringA}
	N(A) = -2 \implies -\frac{4V_{dd}}{\pi A} = -2 \implies A = \frac{2V_{dd}}{\pi} \approx 0.6366 \cdot V_{dd}.
\end{equation}

When $V_{dd} = 1$, solution in (\ref{eq:ringA}) is $0.6366$.
It is reasonably close to the analytical solution $0.618$ in (\ref{eq:ringAstar}).

To acquire the frequency of oscillation, we analyze the linear part of the system.
At the intersection in \figref{ring_nyquist}, 
\begin{eqnarray}
	&& G_2(j\omega)\cdot G_1(j\omega) = -0.5 \implies |\frac{1}{RC\cdot j\omega + 1}| = 0.5 \\
	&& \implies RC\cdot |\omega| = \sqrt{3} \implies |\omega| = \frac{\sqrt{3}}{RC} \\
	&& \implies T = \frac{2\pi}{\sqrt{3}} RC \approx 3.628 \tau \label{eq:ringT}
\end{eqnarray}

Note that result in (\ref{eq:ringT}) is not so close to the analytical solution $T^*$ in
(\ref{eq:ringTstar}).
This is mostly because in this example with ideal relay inverters, the oscillating waveforms
are not close to being sinusoidal, so the results from describing function
analysis are not guaranteed to be accurate.
Nevertheless, an estimation can be acquired for the amplitude and frequency of
the oscillation that are reasonably close the analytical solutions.

In most cases, analytical solutions are not available, but describing function
analysis is a general and systematic approach that can still be applied.
Here, we consider changing the inverter to have more common and realistic
nonlinear property such as
\begin{equation} \label{eq:ringf}
	v_{out} = f(v_{in}) = - \hat A \tanh (k \cdot v_{in}).
\end{equation}

In this case it becomes difficult for analytical solutions to be derived.
But the describing function of (\ref{eq:ringf}) can be easily calculated
analytically or numerically.
Just as in Section \ref{sec:LC}, $N(A)$ can be approximated as
\begin{equation}\label{eq:ringN2}
	N(A) = - \hat A k + \frac{\hat A k^3}{4} \cdot A^2 - \frac{\hat A k^5}{12} \cdot A^4 + o(A^4)
\end{equation}

The block diagram of this ring oscillator is still unchanged from \figref{ringosc3}.
Nyquist plot of $G_2(S)\cdot G_1(s)$ overlaid with $1/N(A)$ is still the same of \figref{ring_nyquist}
except that $N(A)$ is now given as (\ref{eq:ringN2}).
With parameters $\hat A = 1, k = 3$, we can calculate amplitude $A$ at the
intersection in \figref{ring_nyquist}.
\begin{equation} \label{eq:ringA23}
	N(A) = -2 \implies -3 + \frac{3^3}{4} \cdot A^2 - \frac{3^5}{12} \cdot A^4 = -2 \implies A \approx \frac{2}{3}.
\end{equation}

Since (\ref{eq:ringN2}) turns out to be an inaccurate approximation to the
describing function with $\tanh()$, we calculate the describing function
numerically by sweeping $A$ just to double check the result in
(\ref{eq:ringA23}).
The inverter function as in (\ref{eq:ringf}) and the reciprocal of its
describing function are plotted in \figref{ring_tanhNA}.
From \figref{ring_NA} we observe that $1/N(A) = -0.5 \implies A \approx 0.51$.
The linear part of this ring oscillator is unchanged from the one in
\figref{ringosc3}, so the solutions acquired from describing function analysis
are $A \approx 0.51,~T \approx 3.628 \tau$.
\begin{figure}[htbp]
	\subfigure[{Plot of $V_{out}-V_{in}$ relationship as in (\ref{eq:ringf}).
	\figlabel{ring_tanh}}]{
	\epsfig{file=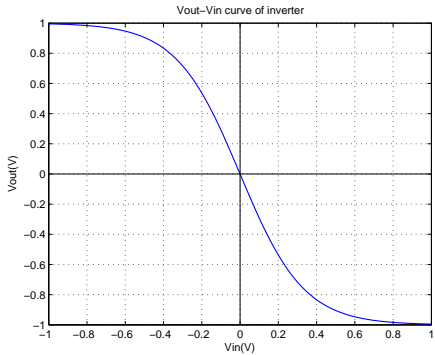,width=0.45\textwidth}
   }\hfill
   \subfigure[{$1/N(A)$ vs. $A$.
	\figlabel{ring_NA}}]{
	\epsfig{file=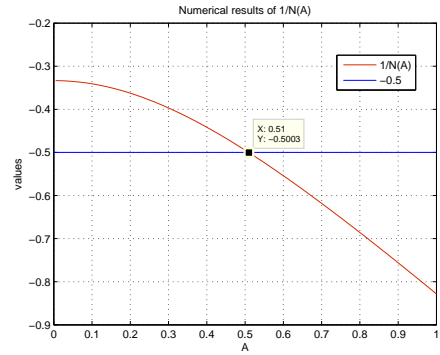,width=0.45\textwidth}
   }
\caption{Plots of characteristics of the inverter as well as the reciprocal of
	its describing function.\figlabel{ring_tanhNA}}
\end{figure}

We can test the results by performing numerical simulation on this 3-stage ring
oscillator, and observe the transient simulation results shown in
\figref{ring_transient}.
Measurements show that the amplitude of the input to the nonlinear block
($e_{in1}$ in \figref{ring_transient}) is $0.50$ and the oscillation period
of the system is around $3.5ms$ (with $\tau = 1ms$ in the circuit).
These properties are quite accurately predicted by the describing function
analysis.
We emphasize again that in this example strict analytical solution is not
available, 
making describing function analysis a valuable tool for providing predictions
for ring oscillators.
\begin{figure}[htbp]
\centering{
	\epsfig{file=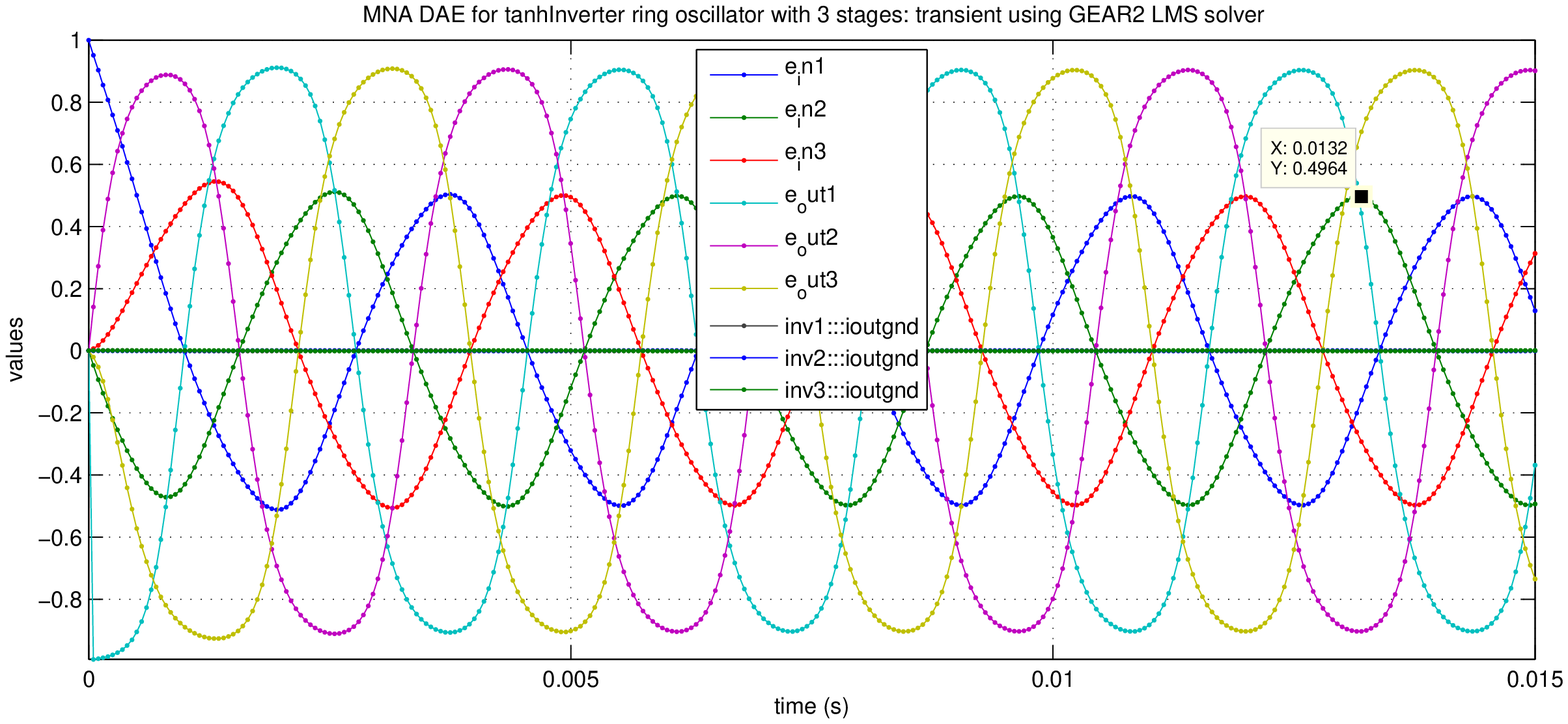,width=\textwidth}
}
\caption{Transient simulation results on system in \figref{ringosc} with
inverter characteristics described in (\ref{eq:ringf}).
    \figlabel{ring_transient}}
\end{figure}

Many oscillators can be categorized as ring oscillators too.
One example is a synthetic biological ring oscillator --- Elowitz
Repressilator. 
Equations for each stage of Elowitz Repressilator can be written as
\begin{eqnarray}
	\frac{d}{dt}m_i &=& -m_i + \frac{\alpha}{1+p_j^n} + \alpha_0 \label{eq:bio_1}\\
	\frac{d}{dt}p_i &=& -\beta (p_i - m_i)\label{eq:bio_2}
\end{eqnarray}

Transient simulation results of Elowitz Repressilator with three stages are
shown in \figref{bio_transient}.
In the simulation, we select parameters as $\alpha = 300$, $\alpha_0 = 0.03$,
$n = 2$, and $\beta = 0.2$.
\begin{figure}[htbp]
\centering{
	\epsfig{file=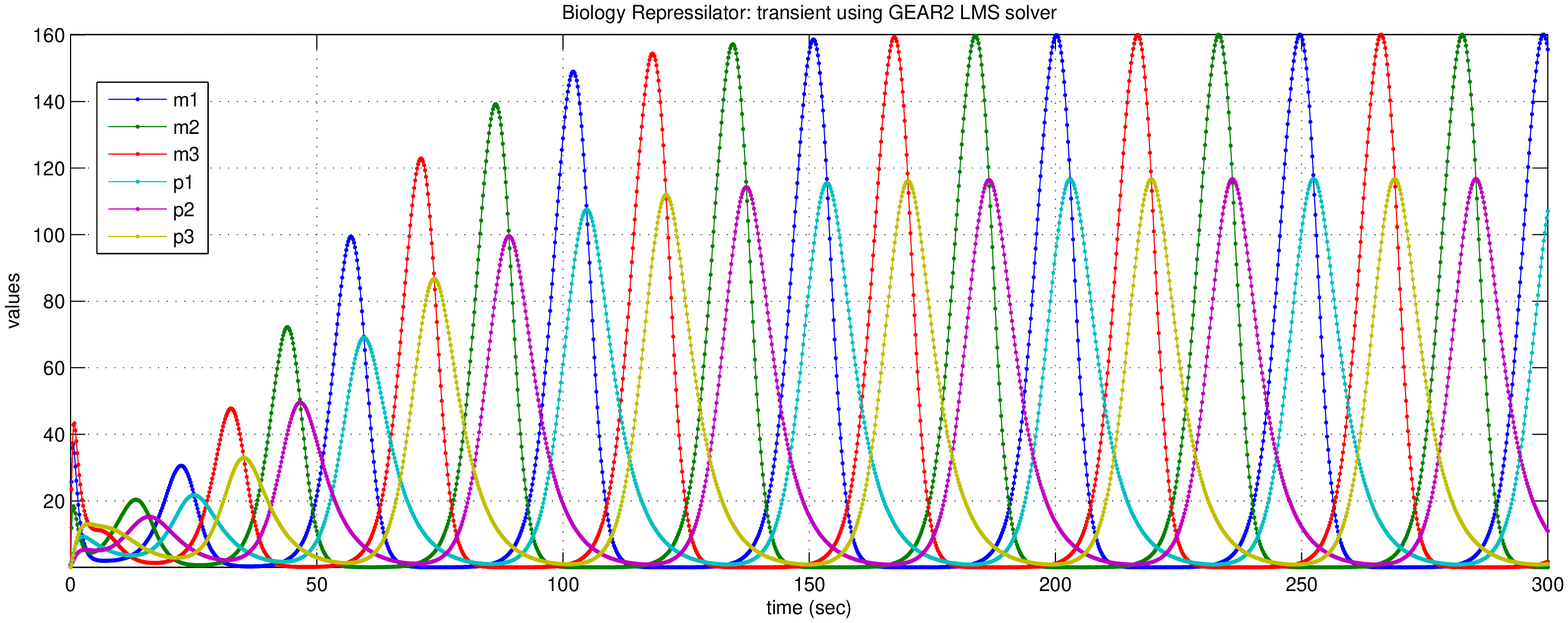,width=\textwidth}
}
\caption{Transient simulation results on Elowitz Repressilator system in (\ref{eq:bio_1}) with (\ref{eq:bio_2})
    \figlabel{bio_transient}}
\end{figure}

Again, this system can be separated into nonlinear and linear blocks.  Block
diagram is the same as \figref{ringosc3} with linear and nonlinear parts:
\begin{equation}
	G_1(s) = \frac{\beta}{s + \beta}, \text{~~}
	G_2(s) = e^{-\frac{2}{3}Ts}
\end{equation}
\begin{equation} \label{eq:bio_f}
	m = f(p) = \frac{\alpha}{1+p^n} + \alpha_0
\end{equation}

Describing function analysis returns approximations to the oscillation
amplitude and frequency --- $T \approx 10/\beta$, $A \approx 40$.
The frequency is rather accurate.
The amplitude indicates that the swing is around 80, while in reality it is
about 115.
Unlike the nonlinearities studied before, function in (\ref{eq:bio_f}) is not
odd symmetric or centered at zero.
So $A\sin(\omega t)$ is not a reasonable approximation to the input of this
nonlinear block any more. 
In order to calculate the describing function of (\ref{eq:bio_f}), offset of
the oscillation has to be estimated.
However, the relationship between this offset and the DC solution of the system
is still not clear to the author.
This limitation of the calculation of describing function requires further
investigation.

\subsection{Relaxation Oscillators} \seclabel{relaxation}

Relaxation oscillation often results from hysteresis.
As is shown in the left plot of \figref{hysteresis}, an electronic relaxation
oscillator is often made by tying RC circuits to a Schmitt trigger in a
feed-back loop.
We can easily separate the nonlinear and linear blocks, then draw the system's
block diagram as in the right plot of \figref{hysteresis}. 
Based on the results in Table \ref{tab:DFs}, the describing function for ideal
hysteresis can be written as
\begin{equation}
    N(A) = \frac{4M}{\pi A} \sqrt{1-(\frac{h}{A})^2} - j \frac{4Mh}{\pi A^2}\cdot, (A \ge h),
\end{equation}
where $M$ is the magnitude of nonlinear function $f(x)$ and $\pm h$ are the transition
points for input $x$.

\begin{figure}[htbp]
\centering{
	\epsfig{file=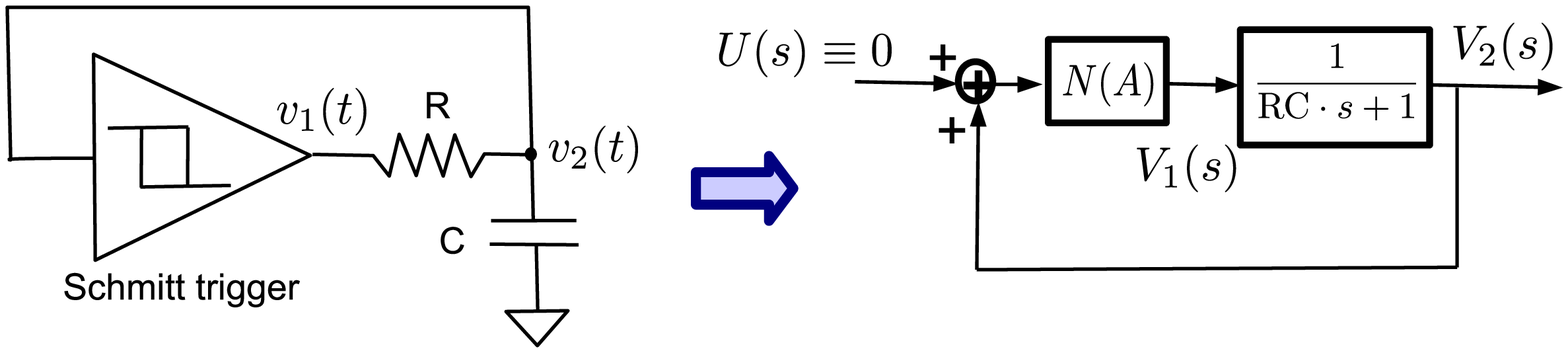,width=0.8\textwidth}
}
\caption{Schematic and block diagram of the relaxation oscillator made of Schmitt
    trigger and RC circuit. \figlabel{hysteresis}}
\end{figure}

The transfer function from $U(s)$ to $V_2(s)$ is 
\begin{equation}
	H(s) = \frac{V_2(s)}{U(s)} = \frac{G(s) \cdot N(A)}{1 - G(s) \cdot N(A)}
\end{equation}

When the system oscillates, $G(s) = 1/N(A)$ has to be satisfied. 
To observe the condition for oscillation more closely, we draw the Nyquist plot
and $1/N(A)$ with different choice of parameters on the complex planes in
\figref{3plots}.
\begin{figure}[htbp]
\centering{
	\epsfig{file=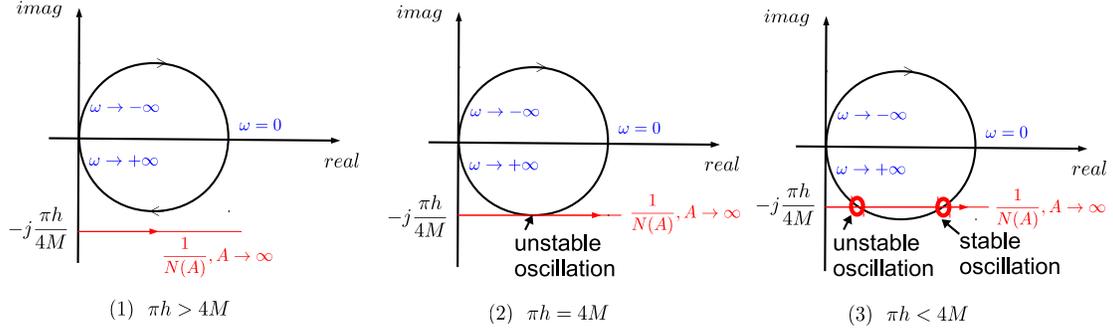,width=\textwidth}
}
\caption{Nyquist plot and the describing function of Schmitt trigger oscillator
plotted on complex plane with conditions $\pi h > 4M$, $\pi h = 4M$ and $\pi h
< 4M$. \figlabel{3plots}}
\end{figure}

Using Nyquist stability criterion as discussed in Section \ref{sec:background}, we analyze
the three scenarios shown in \figref{3plots}.
When $\pi h > 4M$ or $\pi h = 4M$, the system won't oscillate.
When $\pi h < 4M$ there is one stable oscillation operating point.
Oscillation frequency and amplitude can also be calculated from the plot in \figref{3plots}
given all the design parameters.

However, properties of such an ideal relaxation oscillator calculated using describing function
are difficult to be directly verified using simulation.
This is because we modeled hysteresis directly as a memoryless nonlinearity while in fact the
nature of hysteresis indicates the inevitable existence of some form of memory, and should
instead be modelled with differential equations \cite{wang2016well,Wa2017ESD}.

To overcome this difficulty and model hysteresis in a more physical manner, we
model a general relaxation oscillator using equations (\ref{eq:relaxation1})
and (\ref{eq:relaxation2}).
\begin{eqnarray}
	\tau_f\frac{d}{dt} v_o(t) &=& f(v_o(t)) - v_i(t) \triangleq \hat v_i(t) - v_i(t) \triangleq e(t)  \label{eq:relaxation1}, \\
	\tau_s \frac{d}{dt} v_i(t) &=& v_o(t) - v_i(t). \label{eq:relaxation2}
\end{eqnarray}

We then directly convert the system described in (\ref{eq:relaxation1}) and (\ref{eq:relaxation2})
into block diagram and add auxiliary input $U(s)$, as in \figref{relaxation_diagram}.
\begin{figure}[htbp]
\centering{
	\epsfig{file=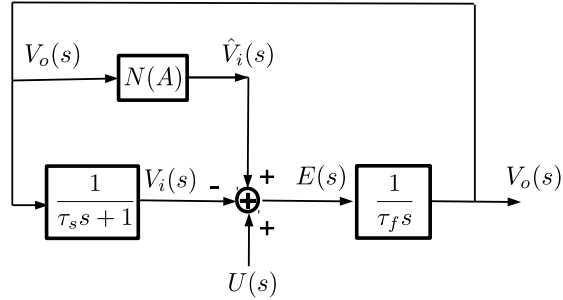,width=0.5\textwidth}
}
\caption{Block diagram of the relaxation oscillator system as described
	in equations (\ref{eq:relaxation1}) and (\ref{eq:relaxation2}). \figlabel{relaxation_diagram}}
\end{figure}

\figref{relaxation_transform} demonstrates the simplification procedures of the system
through transformations on the block diagram in \figref{relaxation_diagram}. 
These simplification procedures are typical in control theory.
After simplification, the nonlinear and linear blocks are separated and we are left with
a simple system that is very similar to the system we studied in Section \ref{sec:background}.
\begin{figure}[htbp]
\centering{
	\epsfig{file=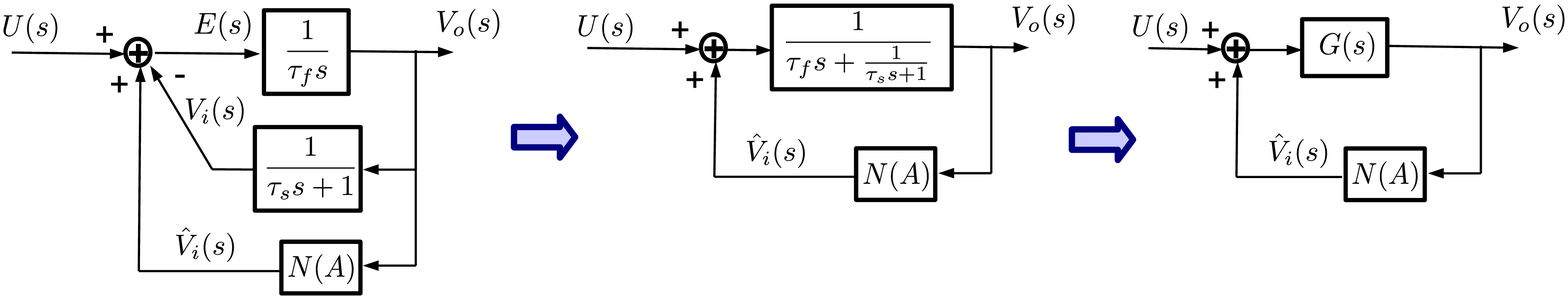,width=\textwidth}
}
\caption{Simplification of block diagram in \figref{relaxation_diagram}.
    \figlabel{relaxation_transform}}
\end{figure}

The transfer function of the resulting system is
\begin{equation}
	H(s) = \frac{V_o(s)}{U(s)} = \frac{G(s)}{1-N(A) \cdot G(s)}
\end{equation}
where
\begin{equation} \label{eq:relaxation_G}
	G(s) = \frac{1}{\tau_f s + \frac{1}{\tau_s s + 1}}.
\end{equation}

For the system to oscillate, $\frac{1}{N(A^*)} = G(j\omega)$ has to be satisfied.
Here we choose our nonlinearity $f$ as
\begin{equation}
	f(v_o) =  -k_1 v_o + k_2 \cdot \tanh(k_3 v_o)
\end{equation}

Analytical calculation result of the describing function of $f$ (approximating $\tanh$ using
its Taylor expansion) is
\begin{eqnarray}
	N(A) &=& \frac{1}{A\pi} \int_0^{2\pi} f(A\sin(\omega t))\cdot \sin(\omega t) d\omega t \\
		 &=& -k_1 + k_2  k_3 - k_2 k_3^3\cdot A^2 + \frac{1}{12}k_2k_3^5 \cdot A^4 + o(A^4).
\end{eqnarray}

In order to analyze oscillation condition quantitatively, we select parameters:
$\tau_f = 2.5 \times 10^{-4}, \tau_s = 1\times 10^{-3}$, 
$k_1=2, k_2=2.5^2, k_3=1/2.5=0.4$.
Then for the system to oscillate, 
$\frac{1}{N(A^*)} = G(j\omega) \implies \frac{1}{N(A^*)} = 4$.
We then plot $\frac{1}{N(A)}$ in \figref{relaxation_N} and observe that
$\frac{1}{N(A^*)} = 4 \implies A^* \approx 1.7$.

\begin{figure}[htbp]
   \subfigure[{Nyquist plot of $G(s)$ as in (\ref{eq:relaxation_G}) overlaid with $1/N(A)$.
	\figlabel{relaxation_G}}]{
	\epsfig{file=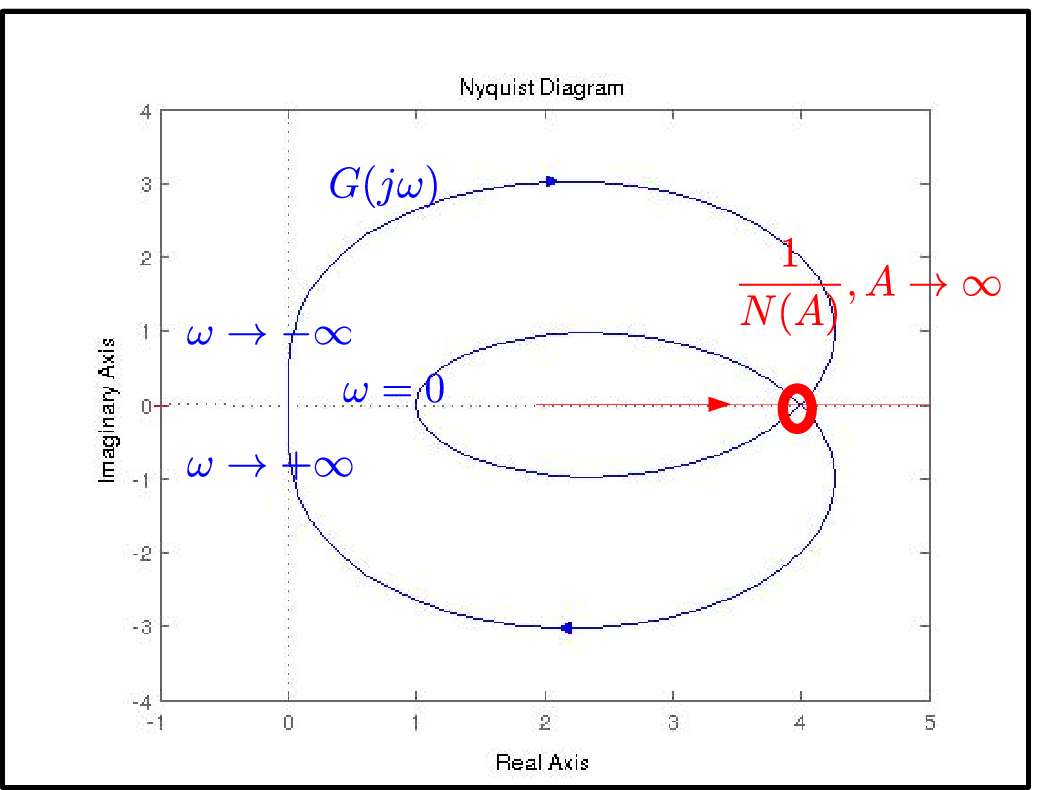,width=0.45\textwidth}
   }\hfill
   \subfigure[{$1/N(A)$ vs. $A$.
	\figlabel{relaxation_N}}]{
	\epsfig{file=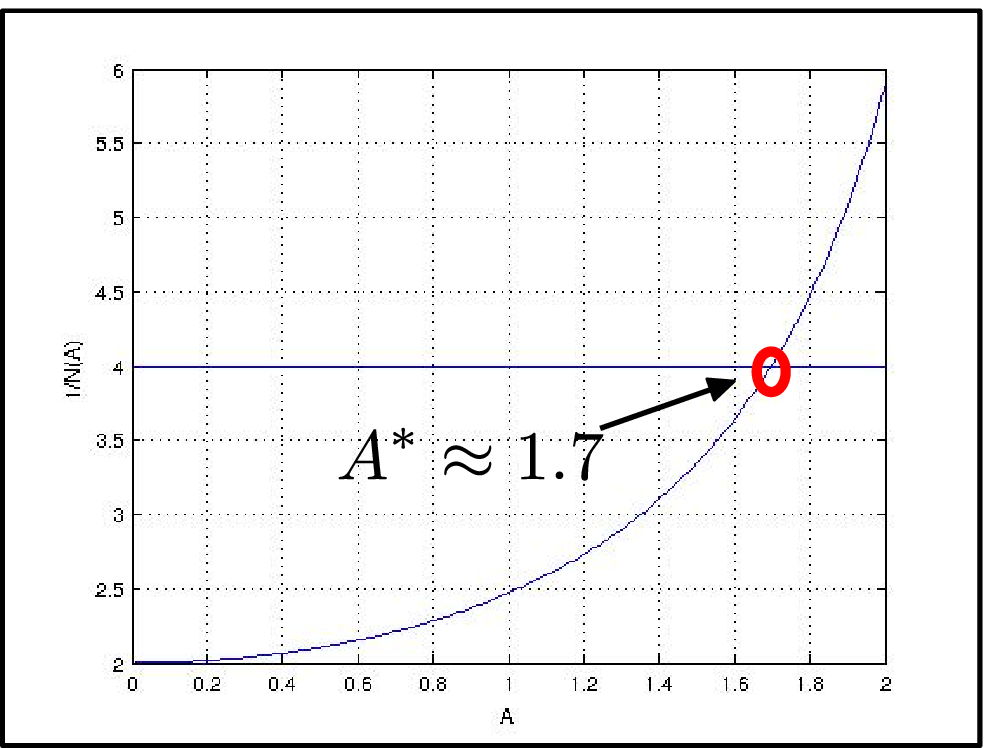,width=0.45\textwidth}
   }
\caption{Nyquist plot and $1/N(A)$ plot of the relaxation oscillator in \figref{relaxation_diagram}.\figlabel{relaxation_GN}}
\end{figure}

The oscillation frequency can be predicted by looking at the intersection in \figref{relaxation_G}.
For $G(j\omega)$ to intersect with real axis, $\angle G(j\omega)$ has to be equal to $0$.
\begin{eqnarray}
	&& \angle G(j\omega)  = \frac{1}{\tau_f \cdot j\omega + \frac{1}{\tau_s \cdot j\omega + 1}} = 0 \\
	\implies && (\tau_s-\tau_f) \omega - \tau_s \tau_f^2 \omega^3 = 0 \\
	\implies && 
    \begin{cases}
        \omega_1 = 0 \\
		\omega_{2,3} = \pm \frac{1}{\tau_s} \sqrt{\frac{\tau_s-\tau_f}{\tau_f}} \label{eq:omega23}\\
    \end{cases}.
\end{eqnarray}

Substitute $\tau_s$ and $\tau_f$ in (\ref{eq:omega23}) with their values, we get the oscillation
frequency $\omega^* \approx 1732$, or $f^* \approx 276 \mathrm{Hz}$.

Then we simulate the system in \figref{relaxation_diagram} and plot the transient waveforms in
\figref{relaxation_transient}.
From the waveforms we can see that the amplitude of the input to the nonlinear block is
approximately $1.65$ while the oscillation frequency is around $260 \mathrm{Hz}$.
These simulation results are close to our predictions acquired through describing function
analysis.

\begin{figure}[htbp]
\centering{
	\epsfig{file=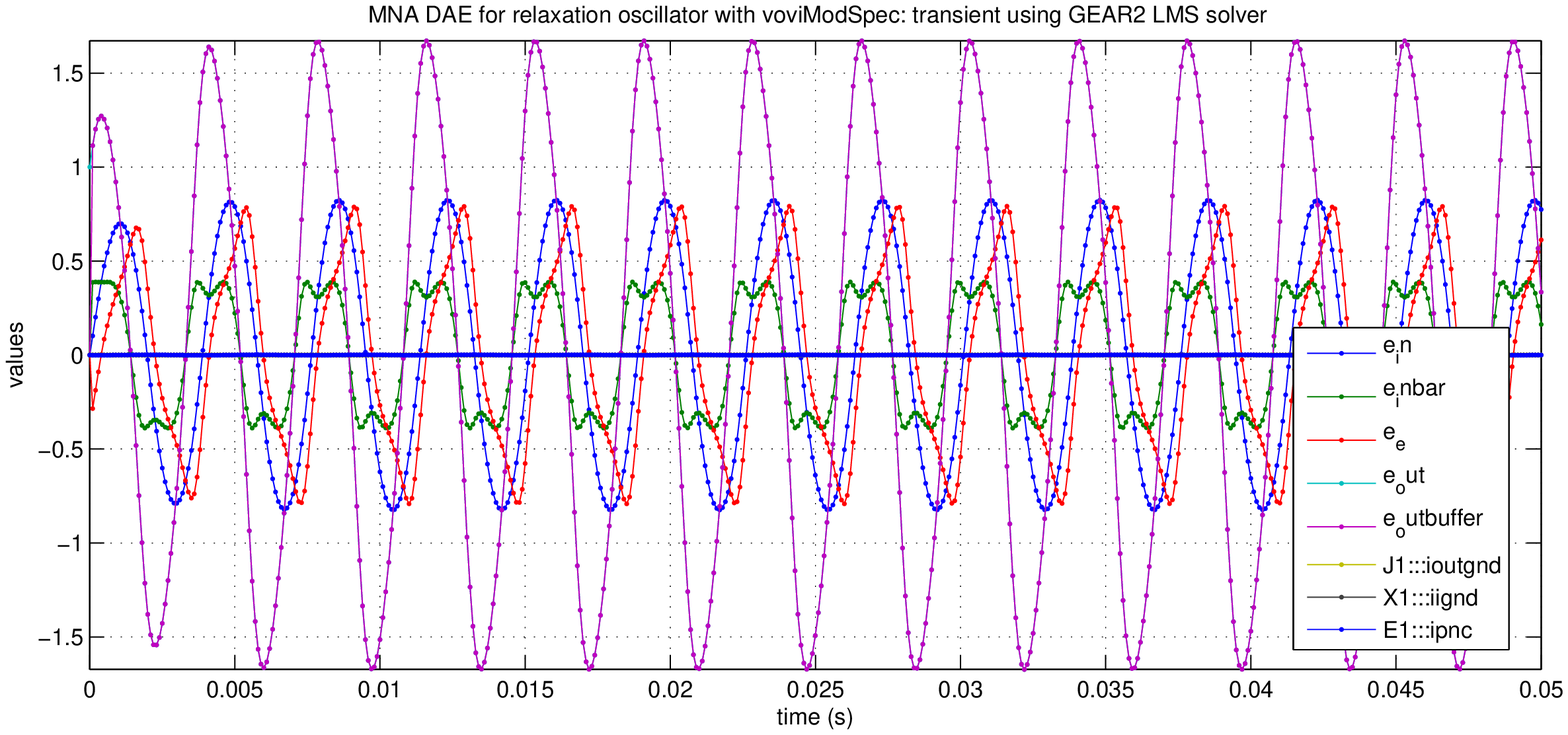,width=\textwidth}
}
\caption{Transient simulation results on system in \figref{relaxation_diagram}.
    \figlabel{relaxation_transient}}
\end{figure}

It is important to note that in the relaxation oscillator model in
(\ref{eq:relaxation1}) and (\ref{eq:relaxation2}), we are not neglecting the
delay from the nonlinear hysteretic component but instead consider both time
constants $\tau_s$ and $\tau_f$.
Through the transformation of block diagrams, we have been able to separate the
linear and nonlinear blocks in general for relaxation oscillators and apply
describing function analysis effectively.

To further demonstrate the use of describing functions, we take a look at
another example of relaxation oscillator from neural biology: Fitzhugh-Nagumo
Neuron Model.
The equations of the system can be written as:
\begin{eqnarray}
	\frac{d}{dt} v &=& v - \frac{v^3}{3} -w + I_{ext} \\
	\tau \frac{d}{dt} w &=& v+a-bw
\end{eqnarray}

\begin{figure}[htbp]
\centering{
	\epsfig{file=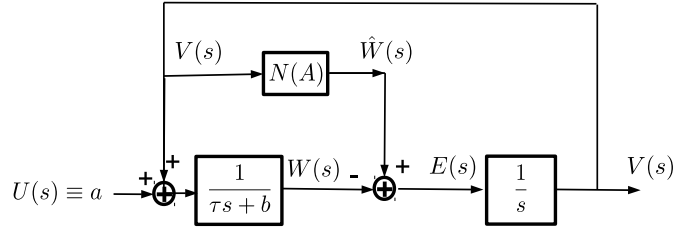,width=0.6\textwidth}
}
\caption{Block diagram of the Fitzhugh-Nagumo neuron relaxation oscillator system.
    \figlabel{relaxation_diagram2}}
\end{figure}

Similar to the Schmitt-trigger-based electronic oscillator, this neuron system can
also be expressed in block diagram format, as is shown in \figref{relaxation_diagram2},
except that $N(A)$ here is the describing function of another $f$:
\begin{equation}
	w =  f(v) = v - \frac{v^3}{3} + I_{ext}
\end{equation}

Its describing function is easy to derive analytically:
\begin{equation}
	N(A) = 1 - \frac{3}{4}A^2.
\end{equation}

Using similar simplification procedures to separate nonlinear and linear blocks and
applying describing function analysis on the system, we will be able to predict 
stable oscillation phenomenon for this neuron system given a set of parameters.
For example, we select parameters $a = 0.7$, $b=0.8$, $\tau=12.5$, $I_{ext} =
0.5$.
Transient simulation shows stable oscillation as in \figref{FNneuron_tran}.
Describing function analysis can be applied in the same procedures as above, 
with $\tau_s = \tau$,  $\tau_f = 1$.
The Nyquist plot is of the same shape as \figref{relaxation_G}, but with the
real axis intersection at $\approx 12.5$ with $\omega \approx 0.25$.
From the intersection we get $A^* \approx 1.1$ and $T^* \approx 25$.
From transient results in \figref{FNneuron_tran}, we observe that the magnitude
of $v$ exceeds 1.1, and the oscillation period is close to 40.
The results from describing function analysis is a little off in this example.

\begin{figure}[htbp]
\centering{
	\epsfig{file=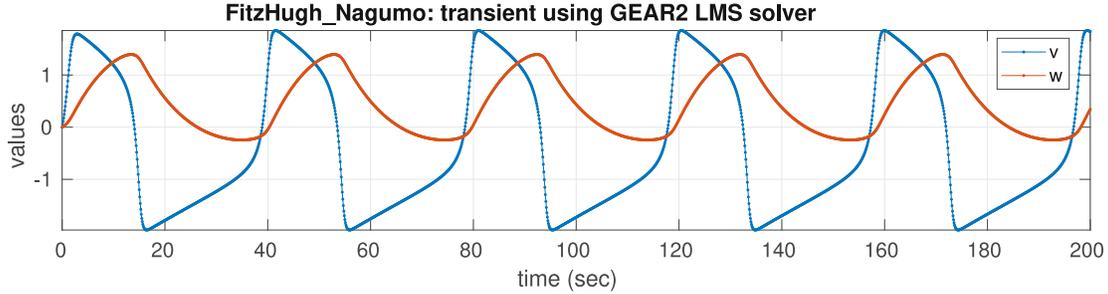,width=\textwidth}
}
\caption{Transient simulation results from Fitzhugh-Nagumo neuron relaxation
oscillator.
    \figlabel{FNneuron_tran}}
\end{figure}

It is worth mentioning that in the case of relaxation oscillators, the system
waveforms are usually not sinusoidal, making the assumptions of describing
function analysis inaccurate.
In simulation, we have noticed that with well-separated $\tau_s$ and $\tau_f$
in the electronic relaxation oscillator system (\ref{eq:relaxation1}) and
(\ref{eq:relaxation2}), the prediction of amplitude and frequency are not
accurate. 
But even in this case, describing function analysis is usually still useful in
predicting the occurrence of oscillation.

\section{Oscillator Design with Describing Functions} \seclabel{design}

Knowledge gained from the analysis of oscillators can often shed light on
the design procedures.
Describing function analysis, with its graphical presentation of systems as
well as its capability to separate systems into linear and nonlinear blocks
with ``orthogonal'' functionalities, often provides useful insights into
oscillator design.
For example, the oscillation frequency involves only the linear block while the
amplitude depends on the nonlinear one; the use of describing function allow
designers to adjust them one by one.
Apart from parameter adjusting, describing function also makes it easier for
people to design oscillators from scratch.
Just like designing negative-resistance LC oscillators, in MEMS 
\cite{miri2010design,nguyen1999integrated,nguyen1995micromechanical}
and spintronics
\cite{dixit2012spintronic,villard2010ghz,finocchio2013nanoscale}
researchers often need to design circuitry around a linear resonator to make
the oscillation self-sustaining.
In these cases, describing function can also serve as a convenient tool.

In this section, we focus on one design example: making relaxation oscillators
behave like a harmonic oscillator, which often indicates better energy
efficiency and less distortion.

\subsection{Design of ``Harmonic'' Relaxation Oscillators} \seclabel{relaxationdesign}

We have studied relaxation oscillators in Section \ref{sec:relaxation}.
Here we look at the oscillator characterized in equations
(\ref{eq:relaxation1}) and (\ref{eq:relaxation2}).

We have already converted the system in to block diagram \figref{relaxation_diagram}.
Now we draw the Bode plot of the linear block in \figref{relaxation_design}.

\begin{figure}[htbp]
\centering{
	\epsfig{file=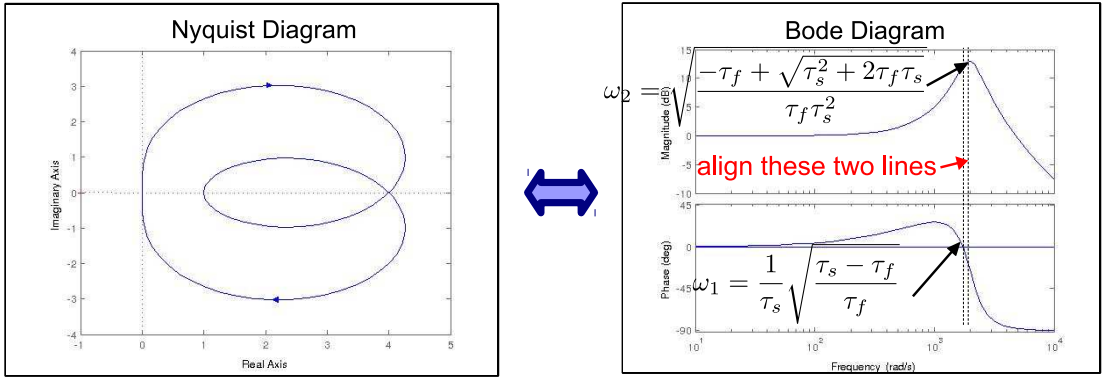,width=\textwidth}
}
\caption{Nyquist and Bode diagrams illustrating design ideas of for the linear block $G(s)$
    to realize ``harmonic-oscillator-like'' relaxation oscillators.
    \figlabel{relaxation_design}}
\end{figure}

Based on \figref{relaxation_G} in Section \ref{sec:relaxation}, oscillation
happens at the location where $\angle G(j\omega_1) = 0$.
Through calculation, we get the expression of $\omega_1$.
\begin{equation}
	\omega_1 = \frac{1}{\tau_s} \sqrt{\frac{\tau_s-\tau_f}{\tau_f}}
\end{equation}

For the waveforms to be more sinusoidal, the linear block should have the
strongest filtering effect around $\omega_1$.
However, based on the Bode plot, the filtering effect is strongest
when the magnitude of transfer function peaks.
Assume it peaks at $\omega_2$, i.e. $\max \{|G(j\omega)|\} = |G(j\omega_2)|$.
\begin{equation}
	\omega_2 = \sqrt{\frac{-\tau_f + \sqrt{\tau_s^2+2\tau_f\tau_s}}{\tau_f\tau_s^2}}
\end{equation}

For $\omega_1$ to be equal to $\omega_2$:
\begin{equation} \label{eq:equate}
	\sqrt{\frac{-\tau_f + \sqrt{\tau_s^2+2\tau_f\tau_s}}{\tau_f\tau_s^2}} = 
	\frac{1}{\tau_s} \sqrt{\frac{\tau_s-\tau_f}{\tau_f}}
	\implies \tau_f = 0 \text{~or~} \tau_s = 0
\end{equation}

(\ref{eq:equate}) indicates under the system structure as in 
(\ref{eq:relaxation1}) and (\ref{eq:relaxation2})
$\omega_1$ will never be equal to $\omega_2$.
So instead, we observe the diagram \figref{relaxation_diagram} and change
the $1/(\tau_s s + 1)$ block to be $1/(\tau_s s + 1)$,
and redraw the diagram in \figref{design_diagram}.
In this way, the transfer function of the linear block becomes
\begin{equation} \label{eq:design_G}
	G(s) = \frac{1}{\tau_f s + \frac{1}{\tau_s s}}
	 = \frac{\tau_s s}{\tau_f\tau_s s + 1}
\end{equation}

Linear system in (\ref{eq:design_G}) has good filtering property so the resulting
waveforms will become more sinusoidal.

It is worth mentioning that the (\ref{eq:design_G}) is essentially the equation of
an LC tank. This indicates mathematically LC oscillators and relaxation oscillators
are very similar.

To implement such a redesigned relaxation oscillator, instead of using RC circuit
that results in $1/(\tau_s s+ 1)$ block, we use integrator that can change the block
transfer function to $1/(\tau_s s)$.
The circuit schematic diagram is given in \figref{design_circuit}.

\begin{figure}[htbp]
\centering{
	\epsfig{file=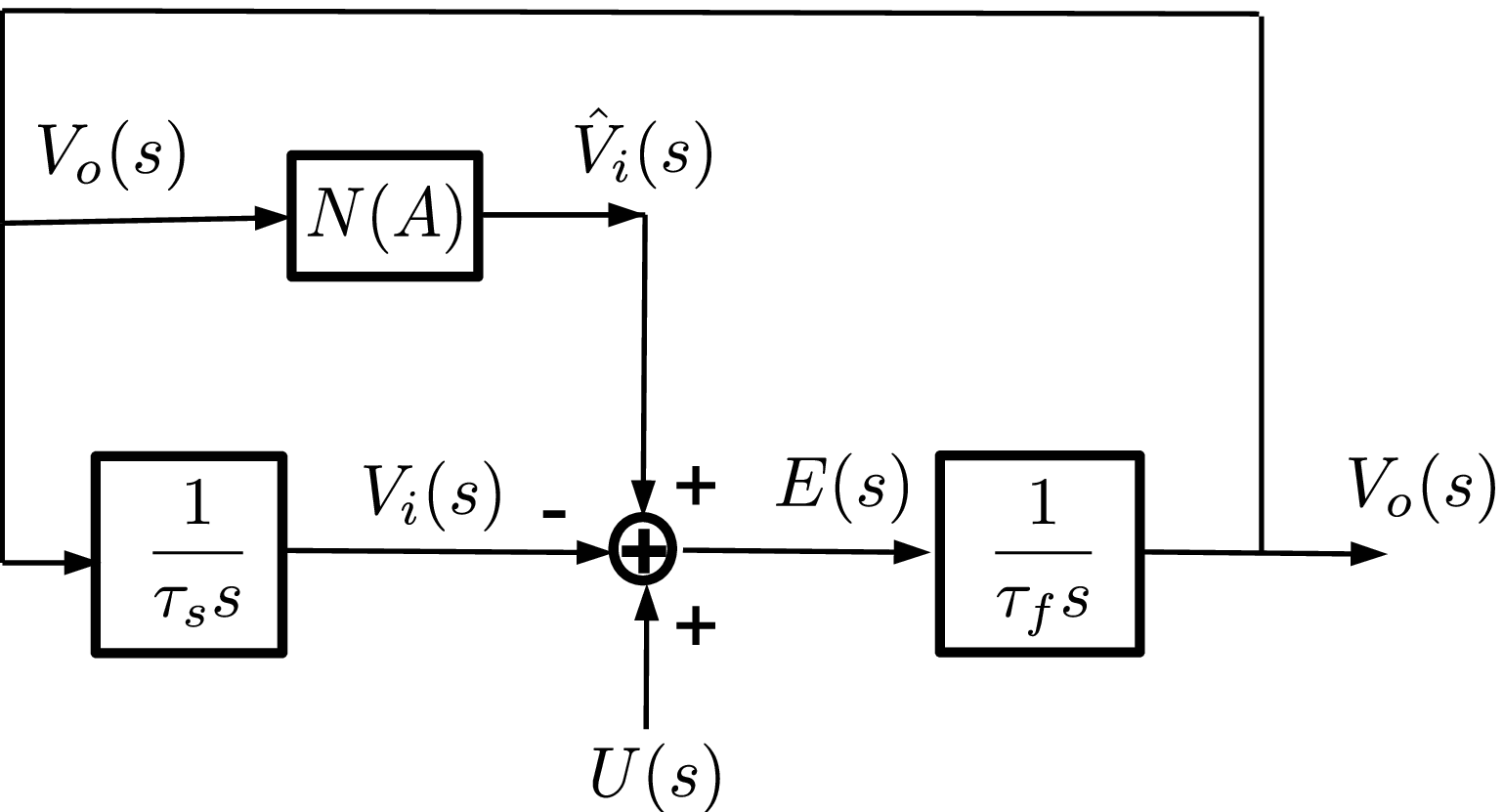,width=0.5\textwidth}
}
\caption{Redesigned diagram of the relaxation oscillator as in \figref{relaxation_diagram}.
    \figlabel{design_diagram}}
\end{figure}

\begin{figure}[htbp]
\centering{
	\epsfig{file=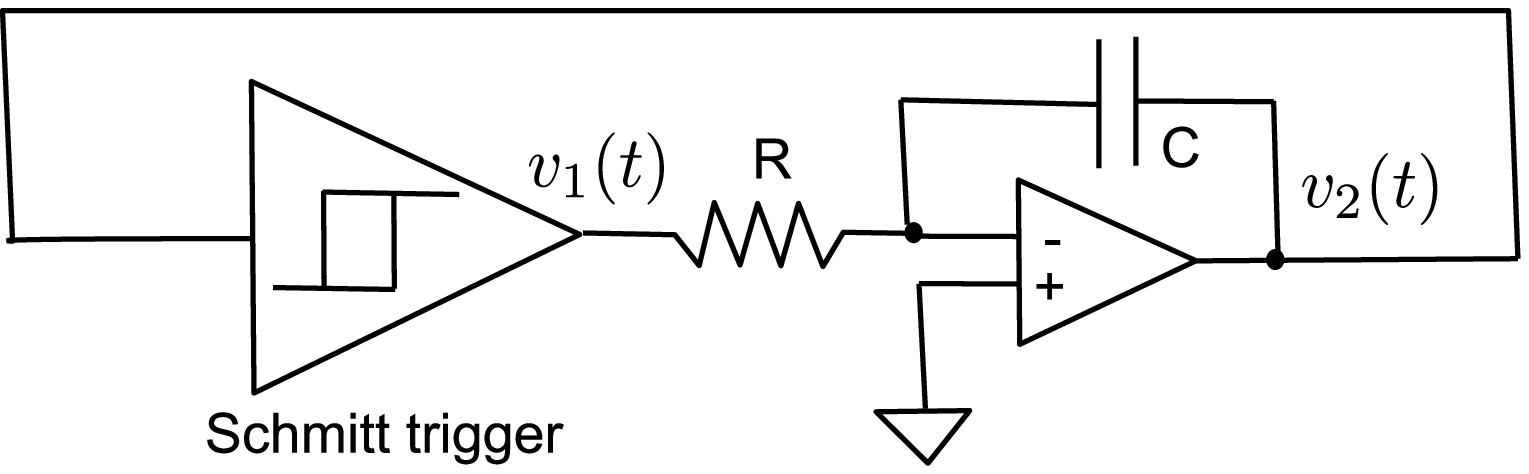,width=0.5\textwidth}
}
\caption{Circuit implementation of the redesigned diagram in \figref{design_diagram}.
    \figlabel{design_circuit}}
\end{figure}

\section{Conclusions}\seclabel{conclusions}

In this manuscript, we studied different types of oscillators using the
describing function analysis.
We showed that the describing function analysis is a general,
domain-independent, systematic approach to the study and design of oscillators.
Examples used in this report demonstrated describing function analysis'
effectiveness, but also revealed its limitations in its support for systems
with non-sinusoidal waveforms.
Such limitations can be overcome by adapting describing functions to
incorporate higher-order harmonics.
Also, describing functions with multiple inputs have potentials in analyzing
multi-tone systems, shedding light on interesting properties such as injection
locking and parametric oscillation.
These are all part of our planned future work on the use of describing
functions in oscillator analysis.

\let\em=\it
\bibliographystyle{unsrt}
\bibliography{stringdefs,jr,PHLOGON-jr,tianshi}

\end{document}